\documentclass[12pt]{amsart}
\usepackage{Baskervaldx, amssymb, graphicx, caption, subcaption, amsmath, array, pdfsync, tikz, tikz-cd, color, url, float, enumerate, geometry, mathrsfs, verbatim, colonequals}
\usepackage[baskervaldx,cmintegrals,bigdelims,vvarbb]{newtxmath}
\usepackage[cal=cm]{mathalfa}
\usepackage[all]{xy}

\usetikzlibrary{calc}

\title{A noncombinatorial proof that toric rank 2 bundles on projective space split}
\author{David Stapleton}

\usepackage[hyperfootnotes=false]{hyperref}
\hypersetup{
  colorlinks= true,
  linkcolor=[RGB]{0,89,42},
  urlcolor=[RGB]{0,89,42},
  citecolor=[RGB]{0,89,42}
}

\theoremstyle{plain}

\newtheorem*{thm}{Theorem}

\theoremstyle{definition}
\newtheorem*{rmk}{Remark}

\newcommand{\hr}[2]{\hyperref[#1]{#2}}

\def\PP{{\mathbf P}}
\def\AA{{\mathbf A}}

\def\CC{{\mathbf C}}

\def\Oc{{\mathcal O}}

\def\PP{{\mathbf{P}}}

\def\Eb{{\mathbb{E}}}

\def\Oc{{\mathcal O}}

\def\Gm{{\mathbb{G}_m}}
\def\Gmn{{(\mathbb{G}_m)^n}}
\def\Gmm{{(\mathbb{G}_m)^m}}

\def\min{{\mathrm{min}}}

\def\ra{{\rightarrow}}

\def\min{{\mathrm{min}}}

\def\Oc{{\mathcal O}}
\def\Ec{{\mathcal E}}

\begin{document}

\maketitle
\thispagestyle{empty}

Let $k$ be an algebraically closed field and let $\Gm$ act on $\PP^n_k$ by scaling the last $n$ coordinates. In this paper we present a new and noncombinatorial proof that $\Gm$-equivariant rank 2 vector bundles on $\PP^n_k$ split when $n\ge 3$.

A conjecture of Hartshorne states that every rank 2 vector bundle $\Ec$ on $\PP^n_k$ splits as a sum of line bundles as soon as $n\ge 7$. We know of two interesting pieces of evidence for this conjecture. First, let $X\subset \PP^n_\CC$ be a smooth complex subvariety which is the zero-locus of a section of $\Ec$. Many of the cohomology groups of $X$ are isomorphic to the cohomology groups of a complete intersection (which is the zero locus of a section of a split bundle). More precisely, Barth proved (\cite{Barth}) that for a nonsingular complex subvariety $X\subset \PP^{n}_\CC$ of codimension $e$, the restriction maps $H^i(\PP^{n}_\CC,\CC)\ra H^i(X,\CC)$ are isomorphisms for $i\le n-2e$. The second piece of evidence comes when one considers toric vector bundles. Considering $\PP=\PP^n_k$ as a toric variety, then a \textit{toric vector bundle} is a vector bundle $\Ec$ on $\PP$ such that the total space $\Eb$ is equipped with a $\Gmn$-action which commutes with the action on $\PP$. Klyachko proved (\cite{Klyachko,Klyachko2}) that every toric vector bundle on $\PP$ of rank $r<n$ splits as a sum of line bundles. (Thus toric rank 2 bundles on $\PP$ split once $n\ge 3$.) Recent work of Ilten and S\"uss (\cite{IltenSuss}) considers smaller tori $\Gmm$ acting on $\PP$ by scaling coordinates. They prove that if $\Ec$ is a rank $r$, $\Gmm$-vector bundle on $\PP$ and $r\le \min\{n,m+3\}$ then $\Ec$ splits as a sum of line bundles. The result of Klyachko (respectively Ilten and S\"uss) is proved by constructing an interesting equivalence of categories between the category of toric vector bundles on a toric variety $X$ and a category of vector spaces with filtration data determined by the fan associated to $X$.

The goal of this note is to give a direct, noncombinatorial proof of the splitting of rank 2 $\Gm$-bundles on $\PP$ for $n\ge 3$. Consider the following scaling $\Gm$-action
\begin{center}
$\Gm\times \PP\ra \PP$\\
\hspace{24pt}$t\times [a_0:a_1:\cdots:a_n]\mapsto [a_0:ta_1:\cdots :ta_n].$
\end{center}
This acts by scaling on the affine space $\AA^n_k=\PP\setminus (x_0=0)$ and fixes the point $[1:0:\cdots:0]$ as well as the hyperplane $(x_0=0)$. In this setting we give a new proof of the following theorem.

\begin{thm}
Let $n\ge 3$. If $\Ec$ is a rank 2, $\Gm$-equivariant vector bundle on $\PP$, then $\Ec$ splits as a direct sum of line bundles.
\end{thm}

\begin{rmk}
As the above action is the restriction to a subgroup of the standard $\Gmn$-action on $\PP$, this recovers Klyachko's result for rank 2 $\Gmn$-equivariant bundles.
\end{rmk}

\begin{proof}
Consider the standard chart $\AA^n_k=(x_0\ne 0)$ and let $H$ denote the hyperplane $(x_0=0)$. Then $\Ec|_{\AA^n}$ is isomorphic to the trivial rank 2 bundle (this does not require the Quillen-Suslin theorem as there is a scaling action on $\Ec|_{\AA^n_k}$). Choose a section $s\in \Ec(\AA^n_k)$ which does not vanish and is an eigenvector for the $\Gm$-action. This gives rise to a meromorphic section of $\Ec$. After possibly twisting by some multiple of $H$ (a $\Gm$-invariant divisor) we can assume $s$ extends to a global section of $\Ec$ (also denoted by $s$) which does not vanish in codimension 1 and which is a $\Gm$-eigenvector. (Note that $\Ec$ is split $\iff$ $\Ec(NH)$ is split for any integer $N$.) Let $Z=(s=0)$ be the scheme-theoretic vanishing locus of $s$.

If we can show $Z$ is a complete intersection we are done. (In this case, $Z$ is the zero locus of a section of $\Oc_{\PP}(a)\oplus \Oc_{\PP}(b)$ for suitable $a$ and $b$. Comparing the two Koszul resolutions for $I_Z$ shows that $\Ec\cong \Oc_{\PP}(a)\oplus \Oc_{\PP}(b)$ as $\mathrm{Ext}^1_{\Oc_{\PP}}(I_Z,\Oc_{\PP}(-a-b))$ is 1-dimensional.)

So we want to understand the ideal $I_Z$. The support of $Z$ is contained in $H$. Consider a point $p\in Z$ and an affine $\Gm$-neighborhood $U$ containing $p$ where $\Ec$ is trivialized. (Such a neighborhood $U$ is an $\AA^1$-bundle over $W=U\cap H$.) Then we have an isomorphism of $\Gm$-algebras:
\[
k[U]\cong k[W][x_0],
\]
where $k[W]$ is trivial as a $\Gm$-representation and $\Gm$ acts on $k\cdot x_0$ with weight 1. As $\Ec$ is trivialized on $U$ we can write $I_Z(U)=(f,g)\subset k[U]$ where $f,g\in k[U]$ are eigenvectors for the $\Gm$-action (i.e. $Z\cap U$ is a complete intersection in $U$). By the above isomorphism of $\Gm$-algebras we have $f=f_W x_0^\alpha$ and $g=g_W x_0^\beta$ where $f_W,g_W\in k[W]$. Without loss of generality (using that $Z$ has codimension 2) we can assume $\alpha=0$. Because $Z$ is supported on $H$, it follows that $g_W$ is a unit in $k[W]/(f_W)$. Thus, we can write:
\[
I_Z(U)=(f_W,x_0^\beta)
\]
for some $f_W\in k[W]$ and $\beta>0$.

We claim that this local description can be extended to global homogeneous equations for $Z$. First, the ideal $(f_W)\subset k[W]$ can be intrinsically defined by $(f_W)=I_Z(U)\cap k[W]\subset k[W]$. Geometrically this means the following: if we consider the map $Z\ra H$ given by linear projection from the point $[1:0:\cdots:0]$, the scheme-theoretic image of the map is an effective Cartier divisor $D\subset H$. Let $F\in H^0(H,\Oc_H(a))$ be the homogeneous equation for $D$. Consider $F$ as a section of $H^0(\PP,\Oc_\PP(a))$. Locally on $U$ we have the equality of ideals $(F|_U) = (f_W)$. Second, it suffices to show that the exponent $\beta$ does not depend on the choice of point $p\in Z$. This is clear for any two points in an irreducible component of $Z$. Because $n\ge 3$, any two irreducible components of $Z$ must intersect (as they are curves in $H$), so $\beta$ is independent of $p\in Z$. Thus $Z$ is a complete intersection defined by the equations $(F=x_0^\beta=0)\subset \PP$.
\end{proof}

\noindent\textit{Acknowledgments.} We would like to thank Morgan Brown, Samir Canning, Kristin DeVleming, Thomas Grubb, and James M\textsuperscript{c}Kernan for helpful conversations.

\bibliographystyle{siam} 
\bibliography{rank2bundles}

\end{document}